\documentclass[12pt,a4paper]{amsart}

\usepackage{amsmath, amsfonts, xifthen, latexsym, amssymb, amsthm, amscd}

\usepackage[utf8]{inputenc}
\usepackage{graphicx}
\usepackage{url}

\usepackage{hyperref}
\hypersetup{colorlinks=true,citecolor=blue,filecolor=blue,linkcolor=blue,urlcolor=blue}
\usepackage[margin=1.4in]{geometry}

\newtheorem{theorem}{Theorem}
\newtheorem{lemma}{Lemma}[section]

\newtheorem{proposition}{Proposition}[section]

\newtheorem{conjecture}{Conjecture}[section]

\newtheorem{corollary}{Corollary}[section]

\theoremstyle{definition}
\newtheorem{definition}{Definition}[section]

\newcommand{\eps}{{\varepsilon}}

\renewcommand{\phi}{{\varphi}}

\newcommand{\wlim}{\operatorname{wlim}}

\newcommand\R{\mathbb R}

\newcommand\Z{\mathbb Z}

\newcommand\opartial{\overline{\partial}}
\newcommand\cG{\mathcal{G}}
\newcommand{\cA}{\mathcal{A}}

\newcommand{\actson}{\curvearrowright}

\newcommand{\cC}{\mathcal{C}}  \newcommand{\cH}{\mathcal{H}}
 \newcommand{\cR}{\mathcal{R}}

\newcommand{\vi}{\vskip 0.1in \noindent}
\usepackage{etoolbox}
\newtoggle{no_cases}\toggletrue{no_cases}
\newcommand{\case}[2][]{\iftoggle{no_cases}{\left\{\begin{array}{ll}#2 & #1}{\\#2 & #1}\togglefalse{no_cases}}
\newcommand{\esac}{\end{array}\right.\toggletrue{no_cases}}

\newcommand{\Prob}{\operatorname{Prob}}

\newcommand{\ut} {\underline{t}}

\newcommand{\Cay}{\operatorname{Cay}}
\newcommand{\Supp}{\operatorname{Supp}}

\begin{document}
\title{Uniform hyperfiniteness}
\author{Gábor Elek}
\address{Department of Mathematics And Statistics, Fylde College, Lancaster University, Lancaster, LA1 4YF, United Kingdom}

\email{g.elek@lancaster.ac.uk}  

\thanks{The author was partially supported
by the ERC Consolidator Grant "Asymptotic invariants of discrete groups, No. 648017 and
by the ERC Starting Grant "Limits of Structures in Algebra and Combinatorics",  No. 805495
}

\begin{abstract}
 Almost forty years ago, Connes, Feldman and Weiss proved that for measurable equivalence relations the notions of amenability and hyperfiniteness coincide. In this paper we define the uniform version of amenability and hyperfiniteness for measurable graphed equivalence relations of bounded vertex degrees and prove that these two notions coincide as well. Roughly speaking, a measured graph $\cG$ is uniformly hyperfinite
 if for any $\eps>0$ there exists   $K\geq 1$ such
 that not only $\cG$, but all of its subgraphs of positive
 measure are $(\eps,K)$-hyperfinite. We also show that
 this condition is equivalent to weighted hyperfiniteness and a strong version
 of fractional hyperfiniteness, a notion recently introduced by
 Lov\'asz.
 As a corollary, we obtain a characterization of exactness of finitely
 generated groups via uniform hyperfiniteness.

\end{abstract}\maketitle
\newpage

\setcounter{tocdepth}{2}
\tableofcontents
\newpage
\section{Introduction}
\noindent
First, let us recall the notion of amenability and hyperfiniteness
in the context of Borel/measurable/continuous combinatorics.\vi
Let $X$ be a standard Borel space. A {\bf Borel graph} $\cG\subset X\times X$ is a Borel set such that
\begin{itemize}
\item for any $x\in X$, $(x,x)\notin \cG$,
\item if $(x,y)\in\cG$, then $(y,x)\in\cG$ as well,
that is, $\cG$ is indeed a graph.(see \cite{KM} for details)

\end{itemize}
 \noindent
 In this paper we always assume that the degrees of a Borel graph
 is countable. The components of a Borel graph $\cG$ are called
 {\bf orbits}. The shortest path metric on the orbits will be
 denoted by $d_\cG$.
 Now, let $\Gamma$ be a countable group
 with symmetric generating system $\Sigma$
 and let $\alpha:\Gamma\actson X$ be a Borel action.
 An associated Borel graph $\alpha^{\Gamma,\Sigma}_G$ is defined 
 in the following way.
 We have $(x,y)\in \alpha^{\Gamma,\Sigma}_G$ for $x\neq y$
 if and only if there is a generator $\sigma\in\Sigma$ such that
 $y\in\alpha(\sigma)(x)$.
 By the Kechris-Solecki-Todorcevic Theorem \cite{KST}, for
 any Borel graph $\cG$, there exists $(\Gamma,\Sigma)$ and an action
 $\alpha:\Gamma\actson X$ such that
 $\alpha^{\Gamma,\Sigma}_G=\cG$. Also, if $\cG$ is
 of bounded vertex degrees, then one can assume that
 $\Gamma$ is finitely generated and $\Sigma$ is a finite generating
 system. \vi
 A {\bf Borel equivalence relation} $E\subset X\times X$ is
 called countable resp. finite, if the equivalence classes are
 countable resp. finite. If $\cG$ is a Borel graph, then the
 associated Borel equivalence relation $E_\cG$ is defined in the
 following way. We have $x\equiv_{E_\cG} y$ if $x$ and $y$ are
 vertices of the same orbit.\vi
 We call the Borel equivalence relation $E$ {\bf hyperfinite} if
 there exist finite Borel equivalence relations $E_1\subset E_2 \subset\dots$ such that $\cup^\infty_{n=1} E_n=E.$ \vi
 We call the Borel equivalence relation $E$ {\bf amenable} if
 there exist Borel functions (the Reiter functions)
 $p_n:E\to [0,1]$ such that
 \begin{itemize}
 \item for any $x\in X$ and $n\geq 1$, $\sum_{z,z\equiv_E x} p_n(x,z)=1,$
 \item for any pair $x\equiv_E y$,
 $$\lim_{n\to\infty} \sum_{z,z\equiv_E x} |p_n(x,z)-p_n(y,z)|=0.$$
 
 \end{itemize}
 \noindent
 It is not hard to see that hyperfiniteness implies amenability.
 However, the converse statement is one of the classical conjectures
 in Borel combinatorics.
 \begin{conjecture}[\cite{KM}]
 Every countable amenable Borel equivalence relation is hyperfinite.
 \end{conjecture}
 \noindent
 Now let us turn to the measurable case (for this part, see also \cite{KM} for further details). Let $(\cG,X)$ be a Borel graph and
 $\mu$ be a Borel probability measure on $X$. We say that $\mu$ is
 an {\bf invariant measure} on $\cG$ if there exists a group
 action $\alpha:\Gamma\actson X$, $\alpha_G^{\Gamma,\Sigma}=\cG$
 such that $\alpha$ preserves the measure $\mu$.
 Note that if $\beta:\Gamma'\actson X$, $\beta_G^{\Gamma',\Sigma'}=\cG$
is another action, then $\beta$ preserves the measure as well.
 Also, if $(\cH,X)$ is another Borel graph such that
 $E_\cH=E_\cG$, then $\mu$ is invariant measure with respect to $\cH$
 as well. 
 \vi
 Similarly, $\mu$ is called a {\bf quasi-invariant} measure on $\cG$
 if the above action $\alpha$ preserves only the measure-class of $\mu$, that is,
 if $\mu(A)=0$ for some Borel set $A\subset X$, then
 $\mu(\alpha(\gamma)(A))=0$ for every $\gamma\in\Gamma$. \vi
 Now, let $(\cG,X,\mu)$ be a {\bf measured graph}, that is,  a Borel graph with a quasi-invariant measure.
 Then, we call $(\cG,X,\mu)$ $\mu$-amenable resp. $\mu$-hyperfinite,
 if there exists a Borel set  
 $Y\subseteq X$ such that
\begin{itemize}
\item $Y$ is a union of equivalence classes of $E_\cG$ (that is, $Y$ is
an {\bf invariant subset}),
\item $\mu(Y)=1$ (that is $Y$ has {\bf full measure}),
\item the induced graph on the set $Y$, $\cG_Y$ is amenable resp. hyperfinite.
\end{itemize}
\noindent
Then we have the celebrated theorem of Connes, Feldman and Weiss \cite{CFW}.
\begin{theorem}
A measured graph $(\cG,X,\mu)$ is $\mu$-amenable if and only if it
is $\mu$-hyperfinite.
\end{theorem}
 
 \noindent Finally, let $(\cG,X,\mu)$ be a bounded degree Borel graph
 with a quasi-invariant measure. For $\eps>0$
 and $K\geq 1$, we call the measured graph $(\cG,X,\mu)$ 
 $(\eps,K)$-hyperfinite if there exists some Borel subset $T\subset X$,
 $\mu(T)<\eps$ such that all the components of the induced
 graph $\cG_{X\backslash T}$ are of size at most $K$.
 Then, $(\cG,X,\mu)$ is $\mu$-hyperfinite if for all $\eps>0$ there
 exists $K\geq 1$ such that $(\cG,X,\mu)$ is $(\eps,K)$-hyperfinite.
 \vi
 Before introducing our new notion, let
 us recall the definition of the {\bf Radon-Nikodym cocycle}.
 Let $(\cG,X,\mu)$ be a measured graph with a quasi-invariant
 measure and $\alpha:\Gamma\actson X$ be a Borel action
 such that $\alpha_G^{\Gamma,\Sigma}=\cG$.
 Then, for any $\gamma\in\Gamma$ we have a Borel function
 $R_\gamma:X\to (0,\infty)$, the Radon-Nikodym derivative, which
 is  unique up to zero-measure perturbation, such that
 for any Borel set $A\subset X$
 \begin{itemize}
 \item $\mu(\alpha(\gamma)(A))=\int_A R_\gamma(x)\,d\mu(x)$,
 \item for any $\gamma,\delta\in\Gamma$ and $x\in X$
 $$R_{\gamma\delta}(x)=R_\gamma(\alpha(\delta))(x) R_\delta(x)\,.$$
 \end{itemize}
 \noindent Hence, we have a Borel function $R:E\to (0,\infty)$,
 the Radon-Nikodym cocycle.
 \begin{definition}
 If for any $\gamma\in\Gamma$ the function $R_\gamma$ is bounded
 and $\cG$ is of bounded vertex degrees, then we call
 $(\cG,X,\mu)$ a measured graph of {\bf bounded type}.
 \end{definition}
 \noindent
 Note the if $\mu$ is an invariant measure, then all the $R_\gamma$'s can be chosen as constant $1$, hence bounded degree graphs with invariant measures are always of bounded type. By the inequality in Section 3.2 of \cite{KV}, for any random walk of a finitely generated group with finitely supported transition measure induces
a bounded type structure on the Furstenberg boundary. 
\noindent
It is important to note that if $\overrightarrow{e}=(x,y)$ is an oriented edge of $\cG$,
we have a well-defined Radon-Nikodym derivative $R_{\overrightarrow{e}}$ corresponding to the edge. If the measured graph is of bounded type, then the
function $R_{\overrightarrow{*}}$ is bounded.
\vi
Before presenting our main definitions, let us
recall the notion of topological amenability for
free continuous actions.
Let $\Gamma$ be a finitely generated
group with symmetric generating set $\Sigma$ and
let $\alpha:\Gamma\actson K$ a free continuous action
of $\Gamma$ on a compact Hausdorff space $K$.
Following \cite{ADR}, we call $\alpha$ \textbf{topologically
amenable} if for any $n\geq 1$, there exists $R_n\geq 1$
and a continuous function $p_n:K\to \Prob(K)$ such that
\begin{itemize}
\item for any $x\in K$, $\Supp(p_n(x))\subset B_{R_n}(x,\alpha_G^{\Gamma,\Sigma})$ (the ball of radius $R_n$
centered at $x$),
\item for all $\alpha_G^{\Gamma,\Sigma}$-adjacent pairs $x,y\in K$
we have
$$\|p_n(x)-p_n(y)\|_1\leq \frac{1}{n}.$$
\end{itemize}
Clearly, the equivalence relation associated to a topologically
amenable action is amenable. Hence, for any quasi-invariant measure 
$\mu$, the measured graph $(\alpha_G^{\Gamma,\Sigma},K,\mu)$
is $\mu$-hyperfinite. By \cite{ADR}, we have the converse if
a the equivalence relation associated to a free continuous action $\alpha$
of a finitely generated group is $\mu$-hyperfinite for
all quasi-invariant measure $\mu$, then $\alpha$ must be
topologically amenable.

\noindent Let $(\cG,X,\mu)$ and $(\cH,Y,\nu)$ be
measured graphs of bounded vertex degrees and
let $\Phi:X\to Y$ be a measure preserving map preserving almost all the orbits. Also, let us assume that there 
exists a constant $L>1$ such that
for $\mu$-almost all $x$ and every $y\equiv_{E_\cG} x$
\begin{equation} \label{coarse}
\frac{1}{L} d_\cG(x,y)< d_\cH(\Phi(x),\Phi(y))<
L d_\cG(x,y).
\end{equation}
\noindent
Then, we say that $\cG$ and $\cH$ are \textbf{coarsely equivalent}.
Note that coarse equivalence is much stronger than
orbit equivalence and $\cG$ is of bounded type if
and only if $\cH$ is of bounded type.
In our paper we introduce a strengthening of the notion
of $\mu$-hyperfiniteness for measured graphs of
bounded type.
\vi
Let $(\cG,X,\mu)$ be a measured graph of bounded type.
Let $A\subset X$ be a Borel subset of positive measure.
Then, one can consider the measured graph
$(\cG_A, A,\mu_A)$ induced on $A$, where for
any Borel set $B\subset A$, 
$$\mu_A(B)=\frac{\mu(B)}{\mu(A)}\,.$$
\noindent
Now we present the key definitions of our paper.
\begin{definition}
Let $(\cG,X,\mu)$ be a measured graph of bounded vertex
degrees, $\eps>0$ and $K\geq 1$. Then, $(\cG,X,\mu)$
is $(\eps,K)$-uniformly hyperfinite if for all subset $A\subset X$ of positive
measure $(\cG_A,A,\mu_A)$ is $(\eps,K)$-hyperfinite.
\end{definition}
\noindent
We call the measured graph $(\cG,X,\mu)$ \textbf{$\mu$-uniformly hyperfinite} if
for any $\eps>0$ there exists $K\geq 1$ such that
$(\cG,X,\mu)$ is $(\eps,K)$-uniformly hyperfinite.
\begin{definition}
Let $(\cG,X,\mu)$ be a measured graph of bounded vertex
degrees, $\eps>0$ and $R\geq 1$. Then, $(\cG,X,\mu)$
is $(\eps,R)$-uniformly amenable, if there exists an invariant
set $Y\subset X$ of full measure, a Borel function $p:Y\to \Prob(Y)$
 such that
\begin{itemize}
\item for all $x\in Y$
$$\Supp(p(x))\subset B_R(x,\cG)\,,$$
\item  and
$$\sum_{x\sim_\cG y} \|p(x)-p(y)\|_1\leq \eps\,.$$
\end{itemize}
\end{definition}
\noindent
We call the measured graph $(\cG,X,\mu)$ of bounded vertex degrees {\bf $\mu$-uniformly amenable}
if for any $\eps>0$ there exists $R\geq 1$ such that
$(\cG,X,\mu)$ is $(\eps,R)$-uniformly hyperfinite.
The main result of the paper is the following theorem.
\begin{theorem}\label{main}
A measured graph $(\cG,X,\mu)$ of bounded type is $\mu$-uniformly \\
amenable if and only if it is $\mu$-uniformly
hyperfinite.

\end{theorem}
\noindent
It will be clear from the proof that for a measured graph
of bounded degrees (without the bounded type condition)
uniform amenability implies uniform hyperfiniteness.
The proof of the theorem will be given by proving
the equivalence of six properties:  uniform amenability,
 uniform local hyperfiniteness, uniform hyperfiniteness,
 weighted hyperfiniteness, approximate strong hyperfiniteness and
 strong fractional hyperfiniteness. 
 \vi
In Section \ref{further}, we give examples of hyperfinite,
but not uniformly hyperfinite measured graphs. Also, 
we show that there exist measured graphs of unbounded type that are 
 uniformly hyperfinite, but not uniformly amenable.
\vi
Finally, in the last section we prove a trichotomy in terms of uniform hyperfiniteness, characterizing
exact non-amenable groups. 

\section{Further motivation and examples} \label{further}
\noindent
It is not very hard to show that $\mu$-hyperfiniteness implies $\mu$-amenability. The converse entails some significant work in the Connes-Feldman-Weiss Theorem. In the case of our Theorem \ref{main}, the more involved part of the proof is to show that $\mu$-uniform hyperfiniteness implies
$\mu$-uniform amenability. In the course of the proof we will show that
$\mu$-hyperfiniteness is equivalent with a series of other notions. 
Let us describe briefly the motivation for this approach. 
Instead of measured graphs, let us consider infinite connected graphs of bounded
vertex degrees. The analogue of $\mu$-uniform amenability is Property A.
This important notion was introduced by Yu \cite{Yu} in the context of the Baum-Connes Conjecture.
\vi
An infinite connected graph $G$ of bounded vertex degrees has \textbf{Property A},
if for any $\eps>0$, there exists $R_\eps>0$ and a function
$p:V(G)\to \Prob(G)$ such that for each $x\in V(G)$,
\begin{itemize}
\item $\Supp(p(x))\subset B_{R_\eps}(x,G)\,,$
\item $\sum_{y,x\sim y} \|p(x)-p(y)\|_1\leq \eps\,.$
\end{itemize}
\noindent
For sequences of graphs hyperfiniteness is well-defined and closely related to $\mu$-hyperfiniteness via the Benjamini-Schramm convergence \cite{Elekamena}. However, there seems to be no sensible way to define  hyperfiniteness for infinite connected graphs. In \cite{Elektimar},
the author and Tim\'ar introduced the notion of weighted hyperfiniteness.
An infinite connected graph $G$ of bounded vertex degrees is weighted
hyperfinite if for any $\eps>0$ there exists $K\geq 1$ such that
for any probability measure $p:V(G)\to [0,1]$, there exists a subset $A\subset V(G)$ such that
\begin{itemize}
\item $p(A)\leq \eps p(V(G))$,
\item the induced graph on $V(G)\backslash A$ has components of size
at most $K$.
\end{itemize}
\noindent
Sako \cite{Sako} proved that Property $A$ is, in fact,
equivalent to weighted hyperfiniteness. Although hyperfiniteness cannot
be defined on a countably infinite graph, one can define
a related notion (this is strongly motivated by the work
in \cite{RWZ}) strong hyperfiniteness. First, recall that
if $G$ is an infinite connected graph $G$ of
bounded vertex degrees, then a subset $A\subset V(G)$ is
a $K$-separator if the induced graph on $V(G)\backslash A$ has components of size
at most $K$. 
\begin{definition}
An infinite connected graph $G$ of bounded vertex degrees is strongly
hyperfinite if for any $\eps>0$ there exists $K\geq 1$ such that we have a probability measure $\nu$ on the compact set
of $K$-separators satisfying the following condition.
For any $v\in V(G)$, the measure of separators containing
$v$ is not greater than $\eps$.
\end{definition}
\noindent
One can prove \cite{Elekuni} that strong hyperfiniteness
is also equivalent to Property A.
The main idea of the proof of Theorem \ref{main} is to show that $\mu$-uniform hyperfiniteness is equivalent to a measured
version of weighted hyperfiniteness, which, in turn, is
equivalent to some measured versions of strong hyperfiniteness and finally, they all are equivalent to $\mu$-uniform amenability.
The steps of the proof are motivated by the proofs of their
combinatorial counterparts.

\noindent
Before starting the proof of our main theorem, let us present
two examples that are intended to demonstrate the subtlety of
the notion of uniform hyperfiniteness.
\vi
\textbf{Example 1.} There exists a measured ergodic graph $(\hat{\cH},\hat{X},\hat{\mu})$ of bounded degrees which is $\hat{\mu}$-hyperfinite,
but not $\hat{\mu}$-uniformly hyperfinite.
\vi
Let $\alpha:\Z\actson X$ the irrational rotation action on the
unit circle preserving the Lebesgue probability measure
$\mu$. We consider the standard generating system
$\Sigma=\{1,-1\}$ and the corresponding measured graph 
$\cH=\alpha_G^{\Z,\Sigma}$. Now, let
$T_1, T_2, T_3,\dots$ be an expander sequence of finite $3$-regular
graphs such that $|V(T_1)|< |V(T_2)| < |V(T_3)|<\dots$.
Let $\{Y_n\}^\infty_{n=1}$
be a sequence of Borel subsets of $X$ satisfying the following
conditions.
\begin{itemize}
\item For any $n\geq 1$, there exists 
an integer $C_n>0$ such that for all $x\in X$ there exists
$y\in Y_n$ so that
$d_{\cH}(x,y)\leq C_n$,
\item $\mu(Y_n) \leq \frac{1}{2^n |V(T_n)|}.$
\end{itemize}
\noindent
The existence of such marker sets is well-known 
(see e.g. \cite{KM}). 
Now, we construct a new measured
graph $(\hat{\cH},\hat{X},\hat{\mu})$
in the following way. First we set
$$\hat{X}=X\cup (Y_1\times V(T_1)) \cup
(Y_2\times V(T_2)) \cup \dots $$
\noindent
Now we define a Borel measure $\nu$
on $\hat{H}$ in the following way.
\begin{itemize}
\item $\nu(A):=\mu(A)$, if $A\subset X$ is
a Borel set.
\item $\nu(B\times \{p\}):=\mu(B)$, if $B\subset Y_n$
is a Borel set and $p\in V(T_n)$ for some $n\geq 1$.
\end{itemize}
\noindent
By our assumption, $\nu(\hat{\cH})\leq 2$.
Now, let $\hat{\mu}$ be the normalized
probability measure associated to $\nu$, that is,
for a Borel set $C\subset \hat{X}$, 
$\hat{\mu}(C):=\frac{\nu(C)}{\nu(\hat{X})}$.
\noindent
Finally, we define a Borel graph
structure on $\hat{X}$. For each $n\geq 1$, fix a vertex
$t_n\in V(T_n)$ and for each $s\in Y_n$, let us connect
$s$ and $s\times t_n$ by an edge $e_s$.
Also, let the induced graph on $s \times V(T_n)$ be $T_n$.
We denote the resulting Borel graph by $\hat{\cH}$, Clearly,
$(\hat{\cH},\hat{X},\hat{\mu})$ is a measured graph with
an invariant probability measure and the 
corresponding orbit equivalence relation is ergodic.
\begin{lemma}
The measured graph $(\hat{\cH},\hat{X},\hat{\mu})$
is $\hat{\mu}$-hyperfinite, but it is not $\hat{\mu}$-uniformly
hyperfinite.
\end{lemma}
\proof
Let $\eps>0$. Since $\cH$ is $\mu$-hyperfinite, we have
a subset $Z\subset X$ and an integer $K\geq 1$ such that
$\mu(Z)< \frac{\eps}{2}$ and all the components of 
$\cH_{X\backslash Z}$ have size not greater than $K$.
\vi
Now, let $q>0$ be an integer
such that
$$Z'= \sum^\infty_{n=q+1} \hat{\mu}(Y_n \times V(T_n)) < \frac{\eps}{2}.$$
\noindent
Then, $\hat{\mu}(Z\cup Z')<\eps$ and 
the size of all the components in the graph $\hat{\cH}_
{\hat{X}\backslash (Z \cup Z')}$ is not greater than
$K+K|V(T_q)|.$ Hence, 
$(\hat{\cH},\hat{X},\hat{\mu})$
is $\hat{\mu}$-hyperfinite. By the expander condition,
for any $l\geq 1$, we have $n_l>0$ such that the graph
$T_{n_l}$ is not $(\eps,n_l)$-hyperfinite. Consequently,
$(\hat{\cH},\hat{X},\hat{\mu})$
is not $\hat{\mu}$-uniformly hyperfinite. So, we have ergodic, invariant
hyperfinite, but non-uniformly hyperfinite measured graphs.
\vi
\textbf{Example 2.} Our second example is
a measured graph $(\hat{H},\hat{X},\hat{\nu})$ which is 
\begin{itemize}
\item $\hat{\nu}$-uniformly
hyperfinite,
\item not $\hat{\nu}$-uniformly amenable,
\item not of bounded type.
\end{itemize}
\noindent
We start with the measured graph 
$(\hat{H},\hat{X},\hat{\mu})$ constructed in Example 1 and substitute
the measure $\hat{\mu}$ with a quasi-invariant
probability measure $\hat{\nu}$ in the same measure class.
Let $T_1, T_2,\dots$ be the graphs in Example 1. and for
$n\geq 1$, let $t_n\in V(T_n)$ be the distinguished vertex.
Also, for $n\geq 1$, let
$t_n=s^n_1,s^n_2,\dots, s^n_{|V(T_n)|}$ be an enumeration
of the vertices of the graph $T_n$.
Finally, we define a probability measure $w_n$ on each vertex set $V(T_n)$
in the following way. \begin{itemize}
\item
Let $w_n(s^n_1)=\frac{1}{2}, w_n(s^n_2)=\frac{1}{4},\dots  $
\item
Let $w_n(s^n_{|V(T_n)|-1})=\frac{1}{2^{|V(T_n)|-1}},
w_n(s^n_{|V(T_n)|})=\frac{1}{2^{|V(T_n)|-1}}.$
\end{itemize}
\noindent
The new measure $\hat{\nu}$ will coincide
with $\hat{\mu}$ on Borel subsets of $X$ and for $n\geq 1$,
$$\hat{\nu}(Y_n\times V(T_n))=\hat{\mu}(Y_n\times V(T_n))\,.$$
We redistribute the weights on each of the sets $(Y_n\times V(T_n))$
in the following way. \vi
Let $n\geq 1$, $1\leq i \leq |V(T_n)|$. Now,
for a Borel set $B_n\subset Y_n$, we define
$$\hat{\nu}(B_n\times \{s^n_i\})=w_n(s^n_i) \hat{\mu}(B_n\times \{s^n_i\}).$$
\noindent
Clearly, $\hat{\mu}$ and $\hat{\nu}$ are in the
same measure class.
So, the measured graph $(\hat{H},\hat{X},\hat{\nu})$ is
still not $\hat{\nu}$-uniformly amenable, since uniform amenability
depends only on the measure class.
Thus, we need to prove the following lemma.
\begin{lemma} \label{examples}
The measured graph $(\hat{H},\hat{X},\hat{\nu})$ is
$\hat{\nu}$-uniformly hyperfinite.
\end{lemma}
\proof 
Let $\eps>0$ and $A\subset \hat{X}$ be a Borel set.
Since, $(\cG,X,\mu)$ is $(\eps,K)$-uniformly hyperfinite,
we have $Y\subseteq A\cap X$ and an integer $K\geq 1$ such that
$\mu(Z)<\eps \mu(A\cap X)$, consequently, 
$\hat{\nu}(Z)<\eps \hat{\nu}(A\cap X)$ and
all the components of $\cG_{(A\cap X)\backslash Z}$
have size at most $K$.
Now, let $j\geq 1$ be an integer such that
\begin{equation} \label{2eq1}
\sum_{i=1}^j \frac{1}{2^i}>1-\eps\,.
\end{equation}
\noindent
Let $Z'$ be the set of elements in $A\backslash X$ in the form
of $y\times t$, where
for some $n\geq 1$, $y\in Y_n$ and $t=s^n_k, k>j$.
Then, by \eqref{2eq1}, $\hat{\nu}(Z) <\eps \hat{\nu} (A\backslash X)$.
Therefore,
$\hat{\nu}(Z\cup Z') <\eps \hat{\nu} (A)$
and the size of the components in $\hat{\cH}_{A\backslash(Z\cup Z')}$
are bounded by $K+Kj$.

\section{Uniform amenability implies uniform local hyperfiniteness}
In the course of our paper we introduce several
notions equivalent to uniform amenability in the realm
of measured graphs of bounded type.
The first such notion is uniform local hyperfiniteness.
Let $(\cG,X,\mu)$ be a measured graph
of bounded vertex degrees. Let $A\subset B\subset X$ be Borel subsets.
Then, the outer boundary set $\opartial_B (A)\subset X$ is defined as the
set of vertices $x\in B\backslash A$ such that there exists $y\in A$, $x\sim_\cG y$.
\begin{definition}
Let $(\cG,X,\mu)$ be a measured graph
of bounded vertex degrees, where $\mu$ is a quasi-invariant
measure. Let $\eps>0, K\geq 1$ be positive constants. We say
that $(\cG,X,\mu)$ is $(\eps,K)$-locally hyperfinite if for any
Borel subset $B\subset X$ of positive measure, there exists
a Borel subset $A\subset B$ of positive
measure such that
\begin{itemize}
\item $\mu(\partial_B (A))< \eps \mu(A),$
\item each component of $\cG_A$ has size at most $K$.
\end{itemize}
\end{definition}
We call the measured graph \textbf{$(\cG,X,\mu)$ $\mu$}-uniformly
locally hyperfinite, if for any $\eps>0$ there
exists $K\geq 1$ such that
$(\cG,X,\mu)$ is $(\eps,K)$-locally hyperfinite.
The main result of this section is the following
proposition.
\begin{proposition}
If the measured graph $(\cG,X,\mu)$ is
$(\eps,R)$-uniformly ame\-nable, then it is
$(\eps,N_{2R})$-locally hyperfinite as well,
where $N_{2R}$ denotes the size of the largest
ball of radius $2R$ in the graph $\cG$.
\end{proposition}
\proof First, we need a lemma.
\begin{lemma}\label{pb}
Let $B\subset Y$ be a Borel set, where $Y$ is the invariant set of
full measure in the definition of $(\eps,R)$-uniform amenability.
Then, there exists a Borel function
$p_B:B\to \Prob(B)$ such that
\begin{itemize}
\item for each $x\in B$, we have
$$p_B(x)\subset B_{2R}(x,\cG)\cap B,$$,
\item for any $x\in B$,
$$\sum_{x\sim_\cG y} \|p_B(x)-p_B(y)\|_1\leq \eps\,.$$
\end{itemize}
\end{lemma}
\proof
If $y\in Y$ is a point that there exists $x\in B$
so that $x\equiv_\cG y$, then we define
$\tau(y)\in B$, $y\equiv_\cG \tau(y)$ such that
$d_\cG(y,\tau(y))=d_\cG(y,B)$
clearly,  the function $\tau$
can be defined in a Borel fashion.
\vi
Now, for $x\in B$ and $z\in B$ we define
the probability measure $p_B(x)$ by setting
$$p_B(x)(z)=\sum_{t, t\in \tau^{-1}(z)} p(x)(t)\,.$$
\noindent
Then $\Supp(p_B(x))\subset B_{2R}(x,\cG)\cap B$
and $\sum_{x\sim y} \| p_B(x)-p_B(y)\|_1\leq \eps$
holds, hence our lemma follows. \qed
\vi
Now we follow the proof of Marks \cite{Marks}. First, we need
a version of Namioka's Trick (Lemma 5.1 \cite{Marks}).
\begin{lemma}
For $a>0$, let $I_{a,\infty}$ be the characteristic 
function of the half-line $(a,\infty)$ and let $f,g\in Prob(X)$
be finitely supported functions. Then,
\begin{equation}\label{831}
\int^\infty_0 \| I_{a,\infty} (f)- I_{a,\infty}(g)\|_1 da= \|f-g\|_1\,.
\end{equation}
\end{lemma}
\noindent
Let $p_B:B\to \Prob(B)$ be the function
defined in Lemma \ref{pb}. Then, we have that
$$\int_B \sum_{x\sim y} 
\int^\infty_0 \| I_{a,\infty} (p_B(x))- I_{a,\infty}(p_B(y))\|_1 \,da d\mu(x)\leq\eps \int_B  \| I_{a,\infty} (p_B(x))\|_1 \,d\mu(x)\,.$$
\noindent
That is, there exists $a\in (0,\infty)$ such that
\begin{equation}\label{832}
\int_B \sum_{x\sim y} 
\| I_{a,\infty} (p_B(x))- I_{a,\infty}(p_B(y))\|_1 \,d\mu(x)\leq\eps \int_B  \| I_{a,\infty} (p_B(x))\|_1 \,d\mu(x)\,.
\end{equation}
\noindent
For $x\in B$, let $\Lambda_x=\{z\,|\, p_B(x)(z)>a\}$. Then, by
\eqref{832},
\begin{equation}\label{MEQ}
\int_B \sum_{x\sim y} |\Lambda_x\triangle \Lambda_y| d\mu(x)\leq\epsilon
\int_B |\Lambda_x| d\mu(x).
\end{equation}
\noindent
By the classical result of Kechris, Solecki and Todorcevic \cite{KST},
there exists a Borel coloring $\phi:X\to Q$ such that
\begin{itemize}
\item $Q$ is a finite set,
\item $\phi(x)\neq \phi(y)$, provided that $d_\cG (x,y)\leq 10 R$.
\end{itemize}
So, for every $x\in B$ and $q\in Q$
there exists at most one $z\in B$
so that
\begin{itemize}
\item $\phi(z)=q$ and
\item either $z\in \Lambda_x$ or $z\in \Lambda_y$ for some
$y$, $x\sim y$.
\end{itemize}
Consequently by \eqref{MEQ}, there exists
an $r\in Q$ such that
\begin{equation}\label{MEQ2}
\int_B \sum_{x\sim y} |\{z \in \Lambda_x\triangle \Lambda_y,\,\phi(z)=r\}|\,d\mu(x)
\leq \eps \int_B |\{z\in \Lambda_x, \phi(z)=r\}|\,d\mu(x).
\end{equation}
\noindent
Let $A$ be the set of the elements
$x\in B$ for which
there exists $z\in B$ such that
$z\in\Lambda_x$ and $\phi(z)=r$.
Observe that the right hand side
of \eqref{MEQ2} equals to $\eps \mu(A)$.
On the other hand, the
left hand side of \eqref{MEQ2} is not
greater than $\mu(\opartial_B(A))$.
Hence, $\mu(\opartial_B(A))\leq \eps \mu(A)$.
Also, all the components of $\cG_A$ has size
at most $N_{2R}$. Thus, our proposition follows. 
\qed
\begin{corollary}\label{cor1}
If the measured graph $(\cG,X,\mu)$ of
bounded vertex degrees is $\mu$-uniformly
amenable, then it is $\mu$-uniformly locally
hyperfinite, as well.
\end{corollary}
\section{Uniform local hyperfiniteness implies uniform hyperfiniteness}

\noindent
The goal of this section is to prove the
following proposition.
\begin{proposition}
$(\eps,K)$-local hyperfinite measured
graphs $(\cG,X,\mu)$ of boun\-ded vertex degrees are $(\eps,K)$-
uniformly hyperfinite.
\end{proposition}
\proof
Let $X_1=X$. By definition, there
exists a Borel set $A_1\subset X_1$ such
that
\begin{itemize}
\item $\mu(\partial_{X_1}(A_1))\leq \eps \mu(A_1)$,
\item all the components of $\cG_{A_1}$ have
size at most $K$.
\end{itemize}
\noindent
Now, let $X_2=X\backslash(A_1\cup \partial_{X_1}(A_1))$ and if let $A_2\subset X_2$ be a Borel set
such that
\begin{itemize}
\item $\mu(\partial_X(A_2))\leq \eps \mu(A_2)$,
\item all the components of $\cG_{A_2}$ have
size at most $K$,
\item $\mu(A_2)>0$, provided that $\mu(X_2)>0.$ 
\end{itemize}
\noindent
By transfinite induction, for each ordinal
we can construct Borel sets $A_\alpha\subset X_\alpha$, such that
\begin{itemize}
\item if $\alpha_1<\alpha_2$ then
$X_{\alpha_1}\supset X_{\alpha_1}$,
\item if $\alpha=\beta+1$, then
$$X_\alpha=X_\beta\backslash (A_\beta\cup 
\partial_{X_\beta}(A_\beta))\,,$$
\item if $\alpha=\lim_{\beta<\alpha}$, then
$X_\alpha=\cap_{\beta<\alpha} X_\beta\,.$
\item $\mu(\partial_X(A_\alpha))\leq \eps \mu(A_\alpha)$,
\item all the components of $\cG_{A_\alpha}$ have
size at most $K$,
\item $\mu(A_\alpha)>0$, provided that $\mu(X_\alpha)>0.$ 
\end{itemize}
\noindent
By our positivity assumption, there exists
a countable ordinal $\alpha$ for which
$\mu(X_\alpha)=0$.
Now, let
$$T:=\cup_{\beta\leq \alpha} \partial_{X_\beta}(A_\beta)\cup X_\alpha$$
\noindent
and let $A=\cup_{\beta\leq\alpha} A_\beta.$
Then,
\begin{itemize}
\item $X\backslash T=A$,
\item $\mu(T)\leq \eps \mu(X)$
\item all the components of $\cG_A$ have
size at most $K$.
\end{itemize}
\noindent
Hence, our proposition follows. \qed
\begin{corollary} \label{cor2}
If the measured graph $(\cG,X,\mu)$ of
bounded vertex degrees is $\mu$-uniformly locally hyperfinite, then
it is $\mu$-uniformly hyperfinite, as well. 
\end{corollary}

\section{Uniform hyperfiniteness implies weighted hyperfiniteness}
\noindent
In Section \ref{further}, we recalled the notion
of weighted hyperfiniteness for countably infinite graphs of bounded vertex degrees. Now,
we introduce a related notion for measured graphs.
\begin{definition}
Let $(\cG,X,\mu)$ be a measured graph of
bounded vertex degrees, $\eps>0$ and $K\geq 1$.
We say that $(\cG,X,\mu)$ is $(\eps,K)$-weighted
hyperfinite, if for any $\mu$-integrable
function $W: X\to [0,\infty)$, there
exists a Borel subset $A\subset X$
such that
\begin{itemize}
\item $\int_A W(x) d\mu(x) \leq
\eps \int_X W(x) d\mu(x)$,
\item all the components of $\cG_{X\backslash A}$
have size at most $K$.
\end{itemize}
\end{definition}
We call $(\cG,X,\mu)$ \textbf{$\mu$-weighted hyperfinite}
if for any $\eps>0$, there exists $K\geq 1$
such that $(\cG,X,\mu)$ is $(\eps,K)$-weighted
hyperfinite.
Clearly, weighted hyperfiniteness implies
uniform hyperfiniteness. The goal of
this section is to prove the converse statement.
Note that the bounded type condition is
crucial.
\begin{proposition}
Let $(\cG,X,\mu)$ be an $(\eps',K)$-uniformly
hyperfinite measured graph of
bounded type, where $d$ is the degree bound of
$\cG$, $M=\sup_{x,y,x\sim y} R_{x,y}$. 
$L=\lceil \frac{3}{\eps} \rceil$, and
$\eps'=\frac{\eps}{3} \left(\frac{3Md}{\eps}\right)^{-L}$.
Then, $(\cG,X,\mu)$ is $(\eps,K)$-weighted
hyperfinite.
\end{proposition}
\proof
We follow the combinatorial approach (used
in the context of finite graphs) by
Romero, Wrochna and \v{Z}ivn\'y \cite{RWZ} up
to the point, where the Radon-Nikodym cocycle
enters the picture.
\vi
So, let $(\cG,X,\mu)$ be an $(\eps',K)$-uniformly
hyperfinite measured graph of bounded type and
$W:  X\to [0,\infty)$ be an integrable
Borel function. Set
$$B_i=\{x\in X\,\mid\, \left( \frac{\eps}{3Md}\right) ^{i+1}
\leq W(x) < \left( \frac{\eps}{3Md}\right) ^{i}\},$$
\noindent
and
for $j\in\{0,1,\dots, L-1\}$
we define
$$B'_j= \cup_{i\in\Z} B_{j+iL}.$$
\noindent
Hence,
we must have $1\leq j^*\leq L-1$ such that
$$W(B'_{j^*})\leq \frac{1}{L} W(X)\leq \frac{\eps}{3}
W(X).$$
\noindent
Now, set
$$C_i= B_{i+j^*+1}\cup B_{i+j^*+2}\cup\dots
\cup B_{i+j^*+L-1}\,.$$
\noindent
Observe that
\begin{equation}\label{e871}
\inf_{x\in C_i} W(x)\geq \left( \frac{\eps}{3Md}\right)^{L} \sup_{x\in C_i} W(x).
\end{equation}
\noindent
Also, if $x\in C_{j}, y\in C_i$ and $i<j$, then
\begin{equation}\label{e872}
W(x)\leq \left( \frac{\eps}{3Md}\right) W(y).
\end{equation}
\noindent
Now, let $F_i\subset X$ be defined as the
set of points $x$ in $X$ such that
$x\in C_j$, $j>i$ and $x\sim y$ for some
$y\in C_i$.
\vi
Then, for any $i\in\Z$, $\mu(F_i)\leq 
Md\mu(C_i)$ and by \eqref{e872},
$$\sup_{x\in F_i} W(x)\leq
\frac{\eps}{3Md}\inf_{y\in C_i} W(y)\,.$$
\noindent
That is,
$$W(F_i)=
\int_{F_i} W(x) d\mu(x)\leq
\frac{\eps}{3} \int_{C_i} W(x) d\mu(x)=\frac{\eps}{3}W(C_i).$$
\vi
Let $F=\cup^\infty_{i=1} F_i$. Then,
$W(F)\leq \frac{\eps}{3} W(X)$.
Now, let $Z=F\cup B'_{j^*}$,
so $W(Z)\leq\frac{2\eps}{3} W(X)$
and consider the graph $\cG_{X\backslash Z}$.
By \eqref{e871}, if $x$ and $y$ are in
the same component of  $\cG_{X\backslash Z}$,
then we have
$$\left(\frac{\eps}{3Md}\right)^L W(y)\leq W(x)\,.$$
\noindent
Since the measured graph $(\cG,X,\mu)$
is $(\eps',K)$-uniform hyperfinite, we
have a set $Z'\subset X\backslash Z$, such that
\begin{itemize}
\item $\mu(Z')\leq\eps'$ and
\item all the components of $\cG_
{X\backslash(Z\cup Z')}$ have size at most $K$.
\end{itemize}
Therefore, we have that
$$W(Z')\leq \eps' \left(\frac{3Md}{\eps}\right)^L
W(X)=\frac{\eps}{3} W(X).$$
Hence, $W(Z\cup Z') \leq \eps W(X)$.
Therefore, the measured graph
$(\cG,X,\mu)$ is $(\eps,K)$-weighted
hyperfinite. \qed
\begin{corollary} \label{cor3}
If the measured graph $(\cG,X,\mu)$ of
bounded type is $\mu$-uniformly hyperfinite, then
it is $\mu$-weighted hyperfinite, as well.
\end{corollary}

\section{Weighted hyperfiniteness implies approximate strong hyperfiniteness}
\noindent
Let $(\cG,X,\mu)$ be a measured
graph of bounded vertex degrees and
$K\geq 1$ be an integer. We say
that $Y\subset X$ is a \textbf{$K$-separator}
if all the components of $\cG_{X\backslash Y}$
have components of size at most $K$.
By the Banach-Alaoglu Theorem the unit ball
$B$ of $L^2(X,\mu)$ is a compact, convex
metrizable space with respect to the
weak topology. If $Y\subset X$, then the
characteristic function of $Y$, $c_Y$ is
an element of $B$. In Section \ref{further}, we
recalled the notion of strong hyperfiniteness
for infinite graphs of bounded vertex degrees, now
we define the obvious analogue of this notion for
measured graphs of bounded vertex degrees.
\begin{definition}
The measured graph $(\cG,X,\mu)$ of bounded vertex degrees is
strongly hyperfinite if
there exists a probability measure $\nu$ on the unit ball
$B$ such that
\begin{itemize}
\item $\nu$ is supported on the characteristic
functions on $K$-separators.
\item The barycenter $b_\nu:=\int_B \,v \,d\nu(v)$ 
satisfies the inequality $b\leq \underline{\eps}$
almost everywhere, where $\underline{\eps}$ is the
constant function taking the value $\eps$.
\end{itemize}
\end{definition}
\noindent
We call $(\cG,X,\mu)$ \textbf{strongly hyperfinite} if for
any $\eps>0$ there exists $K\geq 1$ such that
$(\cG,X,\mu)$ is $(\eps,K)$-strongly hyperfinite.
Unfortunately, we are not able to prove that
strong hyperfiniteness is equivalent to $\mu$-uniform
hyperfiniteness, due to the fact that the set of
characteristic functions of $K$-separators is not closed
(the statement might not even be true). In order to circumvent this difficulty, we introduce two very similar notions which
are, in fact, equivalent to $\mu$-uniform hyperfiniteness:
approximate strong hyperfiniteness and strong fractional 
hyperfiniteness.
\begin{definition}
Let $(\cG,X,\mu)$ be a measured graph of bounded vertex degrees $\eps>0$ and $K\geq 1$.
We say that $(\cG,X,\mu)$ is $(\eps,K)$-approximately strongly
hyperfinite, if there exists a sequence of finitely
supported probability measures 
$$\{\nu_i=\sum_{i=1}^{t_n} x^n_i 
\delta_{c_{Y^n_i}}\}^\infty_{n=1}\,$$
\noindent
non-negative bounded measurable functions $\{z^n\}^\infty_{n=1}$ such that 
$$\wlim_{n\to \infty} \sum_{i=1}^{t_n} (x^n_i c_{Y^n_i}+z^n)=
\underline{\eps}\,,$$
\noindent
where $\wlim$ stands for the weak limit.
\end{definition}
\noindent
Again,  $(\cG,X,\mu)$ is called \textbf{$\mu$-approximately strongly
hyperfinite} if for any $\eps>0$ there exists $K\geq 1$ such that
$(\cG,X,\mu)$ is approximately $(\eps,K)$-strongly hyperfinite.
So, finally we can state the main result of this section.
\begin{proposition}
A measured graph $(\cG,X,\mu)$ of bounded vertex degrees is 
$(\eps,K)$-approximately strongly
hyperfinite  if it is $(\eps,K)$-weighted hyperfinite.
\end{proposition}
\proof
We closely follow the combinatorial proof of Lemma 4.1 in 
\cite{Elekuni}.
Let $C$ be the set of elements $y\in L^2(X,\mu)$
which can be written in the form
$$y=\sum_{i=1}^n t_i c_{Y_i}+z\,,$$
\noindent
where for all $i\geq 1$ $t_i\geq 0, \sum_{i=1}^n t_i=1$
and $z$ is a non-negative function.
\vi
The closure of $C$, $\overline{C}$ is a closed, convex
set in $L^2(X,\mu)$. We have two cases.
\vi
\textbf{Case 1.} $\underline{\eps}\in \overline{C}$.
Then, there exists a sequence of finitely
supported measures 
$$\{\nu_i=\sum_{i=1}^{t_n} x^n_i 
\delta_{c_{Y^n_i}}\}^\infty_{n=1}\,$$
\noindent
together with non-negative bounded measurable functions $\{z_n\}^\infty_{n=1}$ such that 
$$\wlim_{n\to \infty} \sum_{i=1}^{t_n} x^n_i c_{Y^n_i}+z^n=
\underline{\eps}\,.$$
\noindent
that is, $(\cG,X,\mu)$ is $(\eps,K)$-approximately
strongly hyperfinite.
\vi
\textbf{Case 2.} $\underline{\eps}\notin \overline{C}$.
Then, by the Hahn-Banach Separation Theorem
there exists a non-negative $W\in L^2(X,\mu)$
such that
\begin{equation} \label{hahn}
\langle W,\underline{\eps} \rangle <
\langle W, c_Y+z \rangle 
\end{equation}
\noindent
holds for all $K$-separators $Y$ and
bounded non-negative functions $z$.
We can also assume that $W$ is non-negative, otherwise choosing an appropriate $z$ \eqref{hahn} would not hold.
Thus, 
$\int_{c_Y} W(x) d\mu(x)>\eps$ holds
for all $K$-separators $Y$, therefore $(\cG,X, \mu)$ is not
$(\eps,K)$-weighted hyperfinite. \qed
\begin{corollary} \label{cor4}
If the measured graph $(\cG,X,\mu)$ of
bounded type is $\mu$-weighted hyperfinite, then
it is $\mu$-approximately strongly hyperfinite,
as well.
\end{corollary}
\section{Fractional partitions}
\vi
Fractional partitions were recently introduced by Lov\'asz
\cite{Lovasz}. The notion will be crucial in our proof of Theorem \ref{main}, since it will provide the right analogue for strong hyperfinitess of infinite graphs
of bounded vertex degrees.
\vi
Let $(\cG,X,\mu)$ be a measured graph
of bounded vertex degrees. For an integer $K\geq 1$,
we call a subset $L\subset X$ a \textbf{$K$-subset}
if $|L|\leq K$ and the induced graph $\cG_L$ is
connected. Following Lov\'asz let us consider the
Borel space $\cR_K$ of all $K$-subsets.
Note that we have a natural Borel measure $\mu_K$
on $\cR_K$. Let $\cA\subset \cR_K$ be a Borel set, then the Borel function
$\Lambda_\cA:X\to \Z$ is defined in the following way.
$$\Lambda_\cA(x):=|\{A\in \cA\,\mid\, x\in A|\,.$$
\noindent
Then, $$\mu_K(\cA):=\int_X \Lambda_\cA(x)\,d\mu(x)\,.$$
\noindent
For a measurable function $\Phi:\cR_K\to\R$ let
$\Phi^*:X\to\R$ be defined by setting
$$\Phi^*(x):=\sum_{A\in\cR_K,x\in A} \Phi(A)\,.$$
\noindent
Clearly, $\Phi^*$ is a measurable function as well.
A measurable function $\Phi:\cR_K\to \{0,1\}$ is called
a $K$-partition if for all $x\in X$, $\Phi^*(x)=1$.
\begin{definition}
A non-negative measurable function $\Phi:X\to \R$ is
a {\bf fractional $K$-partition} if for almost all $x\in X$, $\Phi^*(x)=1$.
\end{definition}
\noindent
Let $Y\subset X$ be a $K$-separator as in the previous section.
Then, the associated $K$-partition $\Phi_Y$ is defined
in the following way.
\begin{itemize}
\item If $y\in Y$, then $\Phi_Y(\overline{y})=1$,
where $\overline{y}\in\cR_K$ is the singleton containing
$y$,
\item $\Phi_Y(A)=1$, if $A$ is a component of $\cG_{X\backslash Y}$,
\item otherwise, $\Phi_Y(A)=0$.
\end{itemize}
\noindent
Let $\underline{t}=\sum_{i=1}^n t_i \delta_{Y_i}$ be
a finitely supported probability distribution on the
set of $K$-separators. Then, $\Phi_{\underline{t}}:=
\sum_{i=1}^n t_i \Phi_{Y_i}$ is a fractional partition
of $X$.
Now, the measurable function
$\partial\Phi:X\to \R$ is defined by setting
$$\partial\Phi(x)=\sum_{A\in\cR_K,\,x\in\partial A} \Phi(A)\,.$$
\noindent
Note that $\partial A$ (the inner
boundary of $A$) denotes the
set of vertices $x\in A$ such
that there exists $y\notin A$, $x\sim y$.
\vi
For a probability distribution $\underline{t}$ as above
an $\eps>0$, the set $S_{\underline{t},\eps}$ is defined
as the set of points $x\in X$ such that
$$\eps\leq \sum_{i,\,x\in Y_i} t_i\,.$$
\noindent
Also, if $\Phi:\cR_K\to \cR$ is a fractional $K$-partition, then $Q_{\Phi,\eps}$ is defined as the set of points
$x\in X$ such that
$$\eps\leq \partial\Phi(x)\,.$$
\noindent
We end this section with two useful technical lemmas.
\begin{lemma}\label{szerdalemma1}
Let $(\cG,X,\mu)$ be a measured graph of bounded
vertex degrees, let $\underline{t}$ is
a probability distribution as above and $\eps>0$. Then
$$\mu(Q_{\Phi_{\underline{t}},\eps})\leq(d+1)M
\mu(S_{\underline{t},\frac{\eps}{d+1}})\,,$$
\noindent
where $d$ is the degree bound
of $\cG$ and $M=\sup_{x\sim y \in X} R_{x,y}$.
\end{lemma}
\proof If $Y$ is a $K$-separator, then $x\in \partial \Phi_Y$
implies that either $x$ or at least one of its neighbours
is the element of $Y$. Hence,
if $\eps\leq \partial\Phi_{\underline{t}}$, then
there exists $y, d_\cG(x,y)\leq 1$
such that $\sum_{i,y\in Y_i} t_i\geq \frac{\eps}{d+1}.$
So, we have a measurable map $Z:Q_{\Phi_{\underline{t}},\eps}\to
S_{\underline{t},\frac{\eps}{d+1}}$
such that
if $x\in Q_{\Phi_{\underline{t},\eps}}$, then $d_\cG(x,Z(x))\leq 1$.
Therefore, we have
$$\mu(Q_{\Phi_{\underline{t}},\eps})\leq(d+1)M
\mu(S_{\underline{t},\frac{\eps}{d+1}})\,,$$
\noindent
so our lemma follows. \qed
\begin{lemma}\label{szerdalemma2}
Suppose that $\{\Phi_n\}^\infty_{n=1}$
are fractional $K$-partitions such that \\
$\wlim_{n\to\infty} \Phi_n=\Psi$ in
the Hilbert space $L^2(\cR_k,\mu_k)$.
Then, $\Psi$ is a fractional $K$-partition as well.
Also,
\begin{equation} \label{szerda1}
\wlim_{n\to\infty} \partial\Phi_n=\partial \Psi.
\end{equation}
\end{lemma}
\proof
The correspondence $\Phi\to\Phi^*$ defines
a bounded linear map from \\
$L^2(\cR_k,\mu_k)$ onto $L^2(X,\mu)$. Hence,
$\wlim _{n\to\infty} \Phi_n=\Psi$ implies that
$$1=\wlim _{n\to\infty} \Phi^*_n=\Psi^*.$$
\noindent
Thus, $\Psi$ is a fractional $K$-partition.
Similarly, the correspondence $\Psi\to \partial \Psi$ defines
a bounded linear map from
$L^2(\cR_k,\mu_k)$ onto $L^2(X,\mu)$,
therefore \eqref{szerda1} holds as well. \qed
\section{Approximate strong hyperfiniteness implies strong fractional hyperfiniteness}
\begin{definition}
Let $(\cG,X,\mu)$ be a measured
graph of bounded type $\eps>0$, $K\geq 1$.
Then,  $(\cG,X,\mu)$ is $(\eps,K)$-strongly
fractionally hyperfinite
if there exists a fractional $K$-partition $\Phi$
such that for almost all $x\in X$, $\eps\leq \partial \Phi(x).$
\end{definition}
\noindent
We call $(\cG,X,\mu)$ \textbf{strongly fractional $\mu$-hyperfinite} if
for every $\eps>0$ there exists $K\geq 1$ such that
$(\cG,X,\mu)$ is $(\eps,K)$-strongly
fractional hyperfinite.
The goal of this section is to prove the following proposition.
\begin{proposition}
Let $(\cG,X,\mu)$ be a measured
graph of bounded vertex degrees. If
$(\cG,X,\mu)$ is $(\eps,K)$-approximately strongly hyperfinite, then
$(\cG,X,\mu)$ is $(3(d+1)\eps,K)$-strongly
fractional hyperfinite.
\end{proposition}
\proof
Let $\{\underline{t}^n\}^\infty_{n=1}$ be finitely supported distributions on the space of $K$-separators such that
$$\wlim_{n\to\infty} \sum_{i=1}^{s_n} t^n_i c_{Y^n_i}\leq \underline{\eps}\,.$$
\noindent
By definition,
$$\lim_{n\to\infty} \mu(S_{\ut^{n},2\eps})=0\,.$$
\noindent
Hence, by Lemma \ref{szerdalemma1},
\begin{equation} \label{ucso3}
\lim_{n\to\infty} \mu(Q_{\Phi_{\ut^n,2 (d+1)\eps}})=0.
\end{equation}
\noindent
Let $\Psi$ be the weak limit
of a subsequence $\{\Phi_{\ut^{n_k}}\}^\infty_{k=1}$.
By Lemma \ref{szerdalemma2},
$$\wlim_{k\to\infty} \partial \Phi_{\ut^{n_k}}=\partial \Psi$$
\noindent
and 
$\Psi$ is a fractional $K$-partition.
Now, let 
$$A:=\{x\,\mid\,\partial \Psi(x)\geq 3(d+1)\eps\}.$$
\noindent 
Then by weak convergence,
$$2 (d+1)\eps\mu(A)\geq\lim_{k\to\infty} \int_A 
\partial \Psi_{\ut^{n_k}}(x)\,d\mu(x)=
\int_A \partial \Psi(x)\,d\mu(x)\geq 3 (d+1)\eps\mu(A)\,.$$
\noindent
Therefore,
$\mu(A)=0$ and thus $(\cG,X,\mu)$
is $(3(d+1)\eps,K)$-strongly fractionally hyperfinite. \qed
\begin{corollary}\label{cor5}
If the measured graph $(\cG,X,\mu)$ of
bounded vertex degrees is 
$\mu$-approximately strongly
hyperfinite, then it is
$\mu$-strongly fractionally
hyperfinite as well.
\end{corollary}

\section{Strong fractional hyperfiniteness implies uniform amenability}
In this section we finish to proof of Theorem \ref{main}.
\begin{proposition}
Let $(\cG,X,\mu)$ be a measured graph of bounded vertex degrees, $\eps>0$, $K\geq 1$. If $(\cG,X,\mu)$ is $(\eps,K)$-strongly fractionally
hyperfinite, then $(\cG,X,\mu)$ is $(2\eps d,K)$-uniformly amenable, as well.
\end{proposition}
\proof Let $\Phi$ be a fractional $K$-partition 
such that for an invariant subset $Y\subset X$ of full measure
$\eps\leq \partial \Phi(x)$ holds provided that $x\in Y$.
For an element $x\in Y$ let $\Theta(x)$ denote the set of all $K$-subsets
containing $x$. If $A\in\Theta(x)$, let $\phi_A(x):=\frac{1}{|A|}c_A$,
that is, $\sum_{y\in A}\phi_A(y)=1.$
Define $p(x)\in \Prob(Y)$ by setting
$$p(x):=\sum_{A\in\Theta(x)} \Phi(A)\phi_A(x)\,.$$
\noindent
Clearly, $\Supp(p(x))\subset B_K(x,\cG)$.
\begin{lemma}
If $x\sim_\cG y$, then $\|p(x)-p(y)\|_1\leq 2\eps$.
\end{lemma}
\proof
Observe that 
$$\|p(x)-p(y)\|_1\leq \sum_{A,x\in A, y\notin A} \Phi(A)+
\sum_{A,x\notin A, y\in A} \Phi(A)\leq
$$
$$ \leq \sum_{A,x\in \partial A} \Phi(A)+
\sum_{A, y\in \partial A } \Phi(A)\leq 2\eps.
$$
\noindent
Therefore,
$$\sum_{x\sim_\cG y} \|p(x)-p(y)\|_1\leq 2\eps d.$$
\noindent
Hence our Proposition follows. \qed
\begin{corollary} \label{cor6}
If the measured graph $(\cG,X,\mu)$ of
bounded vertex degrees is 
$\mu$-strongly fractionally
hyperfinite, then it is
$\mu$-uniformly amenable, as well.
\end{corollary}
\vi
Now Theorem \ref{main} follows from Corollaries
\ref{cor1}, \ref{cor2}, \ref{cor3}, \ref{cor4},
\ref{cor5} and \ref{cor6}. \qed
\section{Free actions}
\vi
In this section we consider measure class preserving
actions of finitely generated groups such that the Radon-Nikodym derivative of any element is bounded. We call these
actions "actions of bounded type". Also, we call a measure class preserving action uniformly hyperfinite, if the associated measured
graph is uniformly hyperfinite. The following proposition
provides a trichotomy for finitely generated groups. 
\begin{proposition}
Let $\Gamma$ be a finitely generated group. 
\begin{itemize}
\item If $\Gamma$ is amenable, then all free $\Gamma$-actions of
bounded type are $\mu$-uniformly hyperfinite,
\item if $\Gamma$ is non-exact then none of the free $\Gamma$-actions
of bounded type are $\mu$-uniformly hyperfinite,
\item if $\Gamma$ is an exact non-amenable group, then some
of the  free $\Gamma$-actions of bounded type
are $\mu$-uniformly hyperfinite, some of them are not.
\end{itemize}
\end{proposition}
\proof
If $\Gamma$ is amenable, then all free actions of $\Gamma$
are $\mu$-uniformly amenable, hence by our Theorem all
free $\Gamma$-actions of bounded type are $\mu$-uniformly hyperfinite. If $\Gamma$ is non-exact, then none of the
free actions of $\Gamma$ are $\mu$-uniformly amenable, hence
none of the free $\Gamma$-actions of bounded type are
$\mu$-uniformly hyperfinite. If $\Gamma$ is a
non-amenable exact group, then some of its measure-preserving
actions are non-hyperfinite. So, we need to prove that
each such $\Gamma$ has at least one $\mu$-uniformly amenable
actions of bounded type.
\begin{lemma}
Let $\Sigma$ be symmetric generating system for $\Gamma$,
then there exists a non-negative function $\rho:\Gamma\to\R$
and a positive integer $M$ such that
\begin{equation} \label{qw1}
\sum_{\gamma\in\Gamma} \rho(\gamma)=1.
\end{equation}
\noindent
for any $\sigma\in\Sigma$,
\begin{equation} \label{qw2}
\frac{\rho(\gamma\sigma)}{\rho(\gamma)}< M\,.
\end{equation}
\end{lemma}
\proof
Let $\Cay(\Gamma,\Sigma)$ be
the right Cayley graph of $\Gamma$ with the usual 
length function $l(\gamma):=d_{\Cay}(e,\gamma)$,
For $r\geq 0$, let $S_r=\{\gamma\,\mid\,l(\gamma)=r\}.$
Pick a constant $\lambda$ such that $|S_r|\leq e^{\lambda r}$.
Let
$$\rho(\gamma)=\frac{e^{-2\lambda l(\gamma)}}
{\sum_{r=0}^\infty |S_r| e^{-2\lambda l(\gamma)}}\,.$$
\noindent
Then, for large enough $M$ both \eqref{qw1}and \eqref{qw2}
hold. \qed
\vi
Let $\alpha:\Gamma\actson\cC$ be a free continuous
topologically amenable action of $\Gamma$ on the Cantor
set (such action exists by definition). Let $\nu$ be 
the standard Cantor measure. 
The quasi-invariant measure $\mu$ is defined in the usual
way. For a measurable set $A\subset \cC$
$$\mu(A):=\sum_{\gamma\in\Gamma} \rho(\gamma)
\nu(\alpha(\gamma)(A))\,.$$
\noindent
Then, for any $\sigma\in\Sigma$ we have
$$\frac{\mu(\alpha(\sigma)(A))}{\mu(A)}<M.$$
Therefore the action is of bounded type.
Consequently, $\alpha:\Gamma\actson (X,\mu)$ is
a $\mu$-uniformly hyperfinite action. \qed

\end{document}